\documentclass[a4paper]{article}

\usepackage[leqno,intlimits]{amsmath}
\usepackage{amssymb}
\usepackage{amsthm}
\usepackage{hyperref}

\usepackage{pst-all}
\usepackage{pstricks-add}

\usepackage{pstricks}
\usepackage{graphicx}
\usepackage{multido}
\usepackage{calc}

\psset{unit=1pt}
\psset{linewidth=0.5pt}

\newtheorem{tm}{Theorem}[section]
\newtheorem{lm}[tm]{Lemma}
\newtheorem{rem}[tm]{Remark}
\newtheorem{cor}[tm]{Corollary}
\newtheorem{pr}[tm]{Proposition}

\numberwithin{equation}{section}

\newcommand*{\la}{\lambda}
\newcommand*{\ga}{\gamma}
\newcommand*{\al}{\alpha}
\newcommand*{\ba}{\begin{aligned}}
\newcommand*{\ea}{\end{aligned}}
\newcommand*{\be}{\begin{equation}}
\newcommand*{\ee}{\end{equation}}
\newcommand*{\Pv}{{\bf P}}
\newcommand*{\Ev}{{\bf E}}
\newcommand*{\reff}{R^\text{eff}}
\newcommand*{\ceff}{C^\text{eff}}
\newcommand*{\Rb}{\mathbb R}
\DeclareMathOperator{\cp}{cap}

\newcommand{\resist}[1]{
 \qdisk(0,0){2}
 \qdisk(60,0){2}
 \psline(0,0)(15,0)(15,5)(45,5)(45,-5)(15,-5)(15,0)
 \psline(45,0)(60,0)
 \rput(30,0){\footnotesize\({#1}\)}
}

\newcommand{\alap}[2]{
 \qdisk(0,0){2}
 \qdisk(100,0){2}
 \psline(0,0)(5,0)(5,5)(35,5)(35,-5)(5,-5)(5,0)
 \psline(35,0)(47,0)
 \psline(50,-5)(50,5)
 \psline(53,0)(65,0)(65,5)(95,5)(95,-5)(65,-5)(65,0)
 \psline(95,0)(100,0)
 \rput(20,0){\footnotesize\({#1}/2\)}
 \rput(80,0){\footnotesize\({#1}/2\)}
 \psarcn{->}(50,0){10}{150}{30}
 \uput{2}[90](50,10){\footnotesize\(\ast{#2}\)}
}

\newcommand{\alapf}[2]{
 \qdisk(0,0){2}
 \qdisk(100,0){2}
 \psline(0,0)(5,0)(5,5)(35,5)(35,-5)(5,-5)(5,0)
 \psline(35,0)(47,0)
 \psline(50,-5)(50,5)
 \psline(53,0)(65,0)(65,5)(95,5)(95,-5)(65,-5)(65,0)
 \psline(95,0)(100,0)
 \rput(20,0){\footnotesize\({#1}/2\)}
 \rput(80,0){\footnotesize\({#1}/2\)}
 \psarcn{<-}(50,0){10}{150}{30}
 \uput{2}[90](50,10){\footnotesize\(\ast{#2}\)}
}

\newcommand{\dalap}[2]{
 \qdisk(0,0){2}
 \qdisk(100,0){2}
 \psline(0,0)(10,0)(10,5)(30,5)(30,-5)(10,-5)(10,0)
 \psline(30,0)(47,0)
 \psline(50,-5)(50,5)
 \psline(53,0)(70,0)(70,5)(90,5)(90,-5)(70,-5)(70,0)
 \psline(90,0)(100,0)
 \rput(20,0){\footnotesize\({#1}/2\)}
 \rput(80,0){\footnotesize\({#1}/2\)}
 \psarcn{->}(50,0){10}{150}{30}
 \uput{2}[90](50,10){\footnotesize\(\ast{#2}\)}
}

\newcommand{\dalapf}[2]{
 \qdisk(0,0){2}
 \qdisk(100,0){2}
 \psline(0,0)(10,0)(10,5)(30,5)(30,-5)(10,-5)(10,0)
 \psline(30,0)(47,0)
 \psline(50,-5)(50,5)
 \psline(53,0)(70,0)(70,5)(90,5)(90,-5)(70,-5)(70,0)
 \psline(90,0)(100,0)
 \rput(20,0){\footnotesize\({#1}/2\)}
 \rput(80,0){\footnotesize\({#1}/2\)}
 \psarcn{<-}(50,0){10}{150}{30}
 \uput{2}[90](50,10){\footnotesize\(\ast{#2}\)}
}

\newcommand{\primer}[2]{
 \qdisk(0,0){2}
 \qdisk(100,0){2}
 \psline(0,0)(15,0)(15,5)(45,5)(45,-5)(15,-5)(15,0)
 \psline(45,0)(67,0)
 \psline(70,-5)(70,5)
 \psline(73,0)(100,0)
 \rput(30,0){\footnotesize\(#1^\text{pr}\)}
 \psarcn{->}(70,0){10}{150}{30}
 \uput{2}[90](70,10){\footnotesize\(\ast{#2}\)}
}

\newcommand{\primerno}[2]{
 \qdisk(0,0){2}
 \qdisk(100,0){2}
 \psline(0,0)(15,0)(15,5)(45,5)(45,-5)(15,-5)(15,0)
 \psline(45,0)(67,0)
 \psline(70,-5)(70,5)
 \psline(73,0)(100,0)
 \rput(30,0){\footnotesize\(#1\)}
 \psarcn{->}(70,0){10}{150}{30}
 \uput{2}[90](70,10){\footnotesize\(\ast{#2}\)}
}

\newcommand{\secunderf}[2]{
 \qdisk(0,0){2}
 \qdisk(100,0){2}
 \psline(0,0)(15,0)(15,5)(45,5)(45,-5)(15,-5)(15,0)
 \psline(45,0)(67,0)
 \psline(70,-5)(70,5)
 \psline(73,0)(100,0)
 \rput(30,0){\footnotesize\(#1^\text{se}\)}
 \psarcn{<-}(70,0){10}{150}{30}
 \uput{2}[90](70,10){\footnotesize\(\ast{#2}\)}
}

\newcommand{\secunder}[2]{
 \qdisk(0,0){2}
 \qdisk(100,0){2}
 \psline(0,0)(27,0)
 \psline(30,-5)(30,5)
 \psline(33,0)(55,0)(55,5)(85,5)(85,-5)(55,-5)(55,0)
 \psline(85,0)(100,0)
 \rput(70,0){\footnotesize\(#1^\text{se}\)}
 \psarcn{->}(30,0){10}{150}{30}
 \uput{2}[90](30,10){\footnotesize\(\ast{#2}\)}
}

\newcommand{\secunderno}[2]{
 \qdisk(0,0){2}
 \qdisk(100,0){2}
 \psline(0,0)(27,0)
 \psline(30,-5)(30,5)
 \psline(33,0)(55,0)(55,5)(85,5)(85,-5)(55,-5)(55,0)
 \psline(85,0)(100,0)
 \rput(70,0){\footnotesize\(#1\)}
 \psarcn{->}(30,0){10}{150}{30}
 \uput{2}[90](30,10){\footnotesize\(\ast{#2}\)}
}

\begin{document}

\title{Electric network for non-reversible Markov chains}

\author{
 M\'arton Bal\'azs
 \thanks{University of Bristol; part of this work was done while the author was affiliated with the Institute of Mathematics, Budapest University of Technology and Economics; the MTA-BME Stochastics Research Group, and the Alfr\'ed R\'enyi Institute of Mathematics. \texttt{m.balazs@bristol.ac.uk}; research partially supported by the Hungarian Scientific Research Fund (OTKA) grants K60708, F67729, K100473, K109684, and the Bolyai Scholarship of the Hungarian Academy of Sciences.}
 \and
 \'Aron Folly
 \thanks{Department of Mathematics, Ludwig-Maximilians-Universit\"at M\"unchen; part of this work was done while the author was affiliated with the Institute of Mathematics, Budapest University of Technology and Economics; \texttt{follyaron@gmail.com}}
}

\date{\today}

\maketitle

\begin{abstract}
 We give an analogy between non-reversible Markov chains and electric networks much in the flavour of the classical reversible results originating from Kakutani, and later Kem\'eny-Snell-Knapp and Kelly. Non-reversibility is made possible by a voltage multiplier -- a new electronic component. We prove that absorption probabilities, escape probabilities, expected number of jumps over edges and commute times can be computed from electrical properties of the network as in the classical case. The central quantity is still the effective resistance, which we do have in our networks despite the fact that individual parts cannot be replaced by a simple resistor. We rewrite a recent non-reversible result of Gaudilli\`ere-Landim about the Dirichlet and Thomson principles into the electrical language. We also give a few tools that can help in reducing and solving the network. The subtlety of our network is, however, that the classical Rayleigh monotonicity is lost.
\end{abstract}

\small{

}

\section{Introduction}

Random walks or, more generally, reversible Markov chains have a strong connection to electric resistor networks. Our knowledge of this analogy started with the work of Kakutani \cite{kakutani_dirichlet}, Doob \cite{doob_discr_pot}, Kem\'eny, Snell and Knapp \cite{kemeny_snell_knapp}, Nash-Williams \cite{nash-williams_rw}. Since then the field became a foundational part of the theory of reversible Markov chains, we refer the readers to Doyle and Snell \cite{doyle_snell}, Telcs \cite{telcs_art}, Lyons and Peres \cite{lyons_peres_pr}, Chandra, Raghavan, Ruzzo, Smolensky and Tiwari \cite{commt-res} as a few references in the huge literature. Among several results, escape probabilities, transience-recurrence problems, commute and mixing times have been successfully investigated with the use of this analogy. Two fundamental tools were the Thomson (or Dirichlet) energy minimum principles, and Rayleigh's monotonicity law. The former say that under given boundary conditions, the physical current (or voltage, resp.) minimises the power losses on the resistors. As a consequence, Rayleigh's monotonicity law states that the effective resistance of the network is increasing in any of its individual resistances.

The resistor is a symmetric component, and this fact has fundamentally restricted applications to the family of reversible Markov chains. Much less is known therefore in the non-reversible case. The Thomson and Dirichlet principles have been established by Doyle \cite{doyle_energy} and Gaudilli\`ere-Landim \cite{gau-lan}, and re-proved in an elementary way by Slowik \cite{slowik_var}. As an application, Gaudilli\`ere and Landim also prove recurrence theorems in some non-reversible systems. These studies use notions like \emph{energy}, \emph{potential} and \emph{conductance}, but a genuine electric network is not featured behind these ideas. 

In this note we build a full electrical framework behind non-reversible Markov chains. The basic idea is to replace the single resistor by a non-symmetric electrical component. This new part consists of traditional resistors and a new voltage-multiplier unit, which we will just call \emph{amplifier} in short. As shown below, this unit is very directly linked to ``how much a jump is non-reversible'' in the Markov chain. In particular the amplifier becomes trivial and the network reduces to the classical resistor circuit if the chain is reversible. Also, reversing a chain w.r.t.\ its stationary distribution will simply have the effect of reversing the amplifiers.

With this new component, many of the classical analogies work out flawlessly. The starting point is, as in the reversible case, to make voltages in the network directly related to absorption probabilities of the Markov chain. The electric current also has a probabilistic interpretation. Our first observation is that despite the fact that individual components are more complicated than a single resistor, relevant networks can be replaced by a single effective resistance between any to vertices (or even subsets on two different constant potentials). We derive that the effective conductance (reciprocal of the effective resistance) is equal to what people call \emph{capacity} in the theory of Markov chains. We show how symmetry properties and other simple observations regarding the capacity and escape probabilities follow from the electrical point of view. The beautiful observation of Chandra, Raghavan, Ruzzo, Smolensky and Tiwari \cite{commt-res} that connects commute times and effective resistance also generalises without problems to the non-reversible setting. We remark here that a nice mapping of states of non-reversible Markov chains to Euclidean space, based on commuting times, had been worked out earlier by Doyle and Steiner \cite{doyle_steiner_comm_geom}.

Problems start when we look at Rayleigh's monotonicity principle. In its simple naive form monotonicity is just not true in our networks -- this we demonstrate with a counterexample. What can possibly come as a replacement is a question for the future. As a possible first step towards answering this, we rewrite the Dirichlet and Thomson principles by Gaudilli\`ere-Landim \cite{gau-lan} and Slowik \cite{slowik_var} into the electrical language. All terms are then assigned an electrical-energetic meaning, but this has not helped our intuition enough to come up with a sensible way of establishing monotonicity.

Finally, we give a bit of further insight to the behaviour of our electric networks by showing how series and parallel substitutions work. In connection with the lack of (naive) monotonicity, it turns out that delta-star transformations, being essential in the theory of resistor networks, cannot hold in general for our case.

\section{An electric part}

We begin with describing the electric component that we can later use in our analogy with irreversible Markov chains. The schematic picture we use is as follows:
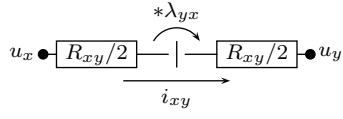
\begin{center}
\begin{pspicture}(120,60)
 \rput(10,30){\alap{R_{xy}}{\la_{yx}}}
 \psline{->}(40,20)(80,20)
 \uput{3}[270](60,20){\footnotesize\(i_{xy}\)}
 \uput{3}[180](10,30){\footnotesize\(u_x\)}
 \uput{3}[0](110,30){\footnotesize\(u_y\)}
\end{pspicture}
\end{center}
This unit is thought of as being connected to neighbouring vertices \(x\) and \(y\) of a graph. These vertices are on respective electric potentials \(u_x\) and \(u_y\), which induces a current \(i_{xy}\) through the unit from vertex \(x\) to vertex \(y\). We will always consider \(i_{\cdot\cdot}\) as an antisymmetric quantity in the sense that \(i_{xy}=-i_{yx}\). Our unit consists of three components: 
\begin{itemize}
 \item an ordinary resistor of \(R_{xy}/2>0\) Ohms,
 \item a \emph{voltage amplifier of parameter \(\la_{yx}>0\)},
 \item another ordinary resistor of \(R_{xy}/2\) Ohms.
\end{itemize}
The resistors each satisfy Ohm's classical law: the signed difference between the potentials at their two ends is proportional to the current that flows through them, and the rate is the value of the resistance. The resistance value \(R_{\cdot\cdot}\) is considered as a symmetric quantity: \(R_{xy}=R_{yx}\). The new element is the \emph{voltage amplifier of parameter \(\la_{yx}\)}. It has the following characteristics:
\begin{itemize}
 \item the current that flows into it on one end agrees with the current that comes out on the other end;
 \item the potential, measured with respect to Ground, from its left end (closer to \(x\)) to its right end (closer to \(y\)) gets multiplied by the positive real parameter \(\la_{yx}\).
\end{itemize}
As the definition naturally suggests, the parameter \(\la_{\cdot\cdot}\) is log-antisymmetric: we always assume \(\la_{yx}=1/\la_{xy}\).  We will also allow the graph to have loops (edges connecting a vertex to itself) from some vertex \(x\) to \(x\), in which case we require \(\la_{xx}=1\).

According to the above, we now follow the potential (with respect to Ground) from left to right in the above unit. First, according to Ohm's law, a drop of \(i_{xy}\cdot R_{xy}/2\) in the potential occurs on the first resistor. Then this dropped potential gets multiplied by \(\la_{yx}\). Finally, a second drop by \(i_{xy}\cdot R_{xy}/2\) occurs on the second resistor. Therefore,
\be
 \ba
  u_y&=\Bigl(u_x-i_{xy}\frac{R_{xy}}2\Bigr)\cdot\la_{yx}-i_{xy}\frac{R_{xy}}2,\qquad\text{or}\\
  i_{xy}&=\frac{2C_{xy}}{1+\la_{yx}}\cdot(\la_{yx}u_x-u_y)
 \ea\label{eq:arram}
\ee
with the introduction of the (Ohmic) conductance \(C_{xy}=C_{yx}=1/R_{xy}\). Notice that the case \(\la_{yx}=1\) reduces our unit to the classical single resistor of value \(R_{xy}\). Notice also that currents are automatically zero along loops: \(i_{xx}=0\) whenever \(x\) is a vertex with a loop.

We write \(z\sim x\) for neighbouring vertices \(z\) and \(x\) in the graph. This includes \(x\sim x\) for vertices \(x\) with a loop. For later use we introduce
\be
 \ga_{xy}:\,=\sqrt{\la_{xy}}=\frac1{\ga_{yx}},\qquad D_{xy}:\,=\frac{2\ga_{xy}}{1+\la_{xy}}C_{xy}=D_{yx},\qquad D_x:\,=\sum_{z\sim x}{D_{xz}\gamma_{zx}}.\label{eq:ddef}
\ee
The symmetry of \(D_{\cdot\cdot}\) follows from that of \(C_{\cdot\cdot}\) and log-antisymmetry of \(\ga_{\cdot\cdot}\) and \(\la_{\cdot\cdot}\). With these quantities we rewrite the above as
\be
 i_{xy}=D_{xy}\cdot(\ga_{yx}u_x-\ga_{xy}u_y).\label{eq:idg}
\ee

We emphasise that the voltage amplifier is not a natural object. Sophisticated engineering would be required to build a black box with this characteristics, and this black box would require an outer energy source (or energy absorber) for its operation. We do not consider this energy source (or absorber) as part of our network.

\subsection{Two alternative parts}

Two alternative units will facilitate calculations in our networks. Using these is not required for any of the later arguments, but simplifies matters. Consider
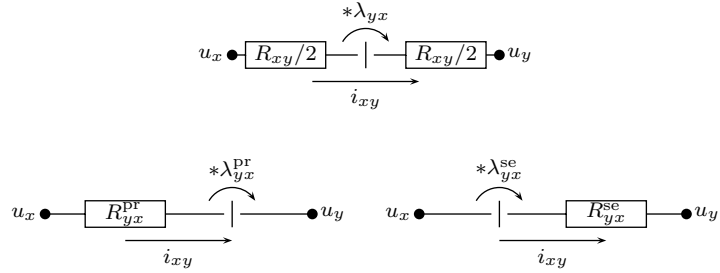
\begin{center}
\begin{pspicture}(260,120)
 \rput(80,90){\alap{R_{xy}}{\la_{yx}}}
 \psline{->}(110,80)(150,80)
 \uput{3}[270](130,80){\footnotesize\(i_{xy}\)}
 \uput{3}[180](80,90){\footnotesize\(u_x\)}
 \uput{3}[0](180,90){\footnotesize\(u_y\)}

 \rput(10,30){\primer{R_{yx}}{\la^\text{pr}_{yx}}}
 \psline{->}(40,20)(80,20)
 \uput{3}[270](60,20){\footnotesize\(i_{xy}\)}
 \uput{3}[180](10,30){\footnotesize\(u_x\)}
 \uput{3}[0](110,30){\footnotesize\(u_y\)}

 \rput(150,30){\secunder{R_{yx}}{\la^\text{se}_{yx}}}
 \psline{->}(180,20)(220,20)
 \uput{3}[270](200,20){\footnotesize\(i_{xy}\)}
 \uput{3}[180](150,30){\footnotesize\(u_x\)}
 \uput{3}[0](250,30){\footnotesize\(u_y\)}
\end{pspicture}
\end{center}
which are the original unit, the \emph{primer} unit and the \emph{secunder} unit, respectively. The primer and secunder units are built of the same types of elements as before. Repeating the arguments we see for the latter two cases
\[
 u_y=(u_x-i_{xy}R^\text{pr}_{yx})\cdot\la^\text{pr}_{yx}\qquad\text{and}\qquad u_y=u_x\la^\text{se}_{yx}-i_{xy}R^\text{se}_{yx}.
\]
Comparing this with \eqref{eq:arram} we conclude that these three units behave in a completely identical way under the choices
\be
 \la_{yx}^\text{pr}=\la_{yx}^\text{se}=\la_{yx},\qquad R^\text{pr}_{yx}=R_{xy}\frac{\la_{yx}+1}{2\la_{yx}},\qquad R^\text{se}_{yx}=R_{xy}\frac{\la_{yx}+1}2\label{eq:prisech}
\ee
which we will assume whenever we write the \(^\text{pr}\) or \(^\text{se}\) indexed quantities. Notice that the primer and secunder resistances are not symmetric quantities anymore.

\subsection{Existence and uniqueness of solutions}
An \emph{electric network} for our purposes consists of our units placed along the edges of a finite, connected graph \(G=(V,\,E)\). We allow \(G\) to have loops as well. Suppose that a subset \(W\) of the vertices is taken to fixed potentials, \(U_x\), \(x\in W\). The only requirement we make is that \(W\) is non-empty. We show below that there exists a unique solution of the network with these boundary values that is, a unique set of currents \(i_{\cdot\cdot}\) with
\be
 \sum_{y\sim x}i_{xy}=0\qquad\forall x\notin W,\label{eq:dfree}
\ee
and voltages \(u_x\) with
\[
 u_x=U_x\qquad\forall x\in W
\]
that satisfy \eqref{eq:arram} for all \(x\sim y\). We start with uniqueness.
\begin{pr}
 Given the graph \(G\) and the boundary set \(W\), fix the boundary condition \((U_x)_{x\in W}\), and suppose that we have two solutions \(u'_\cdot\), \(u_\cdot\) and \(i'_{\cdot\cdot}\), \(i_{\cdot\cdot}\) with this boundary condition. Then \(u'\equiv u\) and \(i'\equiv i\).
\end{pr}
\begin{proof}
 As \eqref{eq:arram} is linear, the difference of two solutions is yet another solution. Therefore \(u'_\cdot-u_\cdot\) and \(i'_{\cdot\cdot}-i_{\cdot\cdot}\) is another solution with boundary condition 0 for all \(x\in W\). Define now the set \(\oplus\subset V-W\) of vertices where \(u'-u\) is positive. If this set is nonempty, then any edge that connects it with the rest of \(V\) sees an outflowing current by \eqref{eq:arram}. But this contradicts \eqref{eq:dfree} (summed up for \(x\in\oplus\)). A similar argument shows that there are no vertices of negative potential either, thus \(u'-u\equiv0\). Then by \eqref{eq:arram} it follows that \(i'-i\equiv0\) as well.
\end{proof}
We proceed by existence of solutions. Call the \emph{incoming current to a vertex \(x\in W\)} \(i_x\):
\be
 i_x:\,=\sum_{y\sim x}i_{xy},\label{eq:ic}
\ee
this is zero for all \(x\notin W\).
\begin{lm}
 Given the graph \(G\) and the boundary set \(W\), suppose we have a solution for all boundary conditions \((U_x)_{x\in W}\). Then for any \(x\in W\), \(i_x\) is an affine increasing function of \(U_x\) when keeping all other boundary voltages \(U_y\), \(x\ne y\in W\) constant.
\end{lm}
\begin{proof}
 Fix \(x\in W\), consider \(U'_x>U_x\), \(U'_y=U_y\) for all \(x\ne y\in W\) and the corresponding solutions \(u'_\cdot\), \(u_\cdot\) and \(i'_{\cdot\cdot}\), \(i_{\cdot\cdot}\). Then \(u'_\cdot-u_\cdot\) and \(i'_{\cdot\cdot}-i_{\cdot\cdot}\) is another solution with boundary condition \(U'_x-U_x>0\) for \(x\), and 0 for all \(x\ne y\in W\). Again looking at the current out of the set \(\oplus\) of vertices with positive potential it is clear that the incoming current \(i'_x-i_x\) is strictly positive in this setting. Now, multiplying all of \(u'-u\), \(i'-i\) by any factor \(\beta\) is yet another solution. It follows that \(i'_x-i_x\) is a positive constant multiple of \(U'_x-U_x\) and the proof is complete.
\end{proof}
\begin{pr}
 There is a solution for any set \(\emptyset\ne W\in V\) and boundary condition \((U_x)_{x\in W}\).
\end{pr}
\begin{proof}
 We will call the set \(V-W\) the free vertices, and perform an induction on its size \(n=|V-W|\). When \(n=0\), all potentials are fixed, and the currents are simply computed by \eqref{eq:arram}. Suppose now that the statement is true for \(n\), and consider a set \(W\) with \(|V-W|=n+1\). Pick any vertex \(x\in |V-W|\). Fixing all boundary values \(U_y\) for \(y\in W\) and also the value \(U_x\), we only have \(n\) free vertices, and we know by the induction hypothesis that we have a solution. We also know by the above lemma that the incoming current \(i_x\) is an affine function of \(U_x\). Therefore there exists a particular value \(U_x^0\) with the corresponding incoming current \(i_x^0=0\), and a solution \(u^0_\cdot\), \(i^0_{\cdot\cdot}\) that goes with the boundary condition \(U^0_x\) for \(x\) and \(U_y\) for \(y\in W\). This will be a solution with boundary condition \(U_y\) for \(y\in W\) only, and the induction step is complete.
\end{proof}

\section{Irreversible Markov chains}

\subsection{Absorption probabilities and the connection}

In this section we make a connection of electric networks, built of our units, to Markov chains. The novelty is that the chain does not need to be reversible. For the following proposition two non-intersecting subsets, \(A\) and \(B\) of the vertex set \(V\) are supposed to be connected to constant external potentials \(U_A\) and \(U_B\):
\[
 u_a\equiv U_A\qquad\forall a\in A\qquad\text{and}\qquad u_b\equiv U_B\qquad\forall b\in B.
\]
All other vertices are free: they just connect to neighbouring ones via our units. The starting point is
\begin{pr}\label{pr:uuu}
 For every \(x\notin A\cup B\), we have
 \be
  u_x=\sum_{y\sim x}u_y\frac{D_{xy}\gamma_{xy}}{D_x}.\label{eq:uuu}
 \ee
\end{pr}
\begin{proof}
 Consider a vertex \(x\notin A\cup B\), and its neighbouring vertices \(y\sim x\), \(x\ne y\). We demonstrate the situation for two neighbours with the picture
 
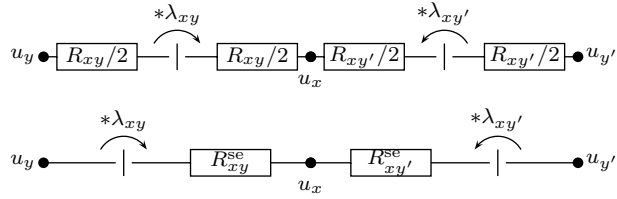
\begin{center}
  \begin{pspicture}(220,80)
   \rput(10,60){\alap{R_{xy}}{\la_{xy}}}
   \rput(110,60){\alapf{R_{xy'}}{\la_{xy'}}}
   \uput{3}[270](110,56){\footnotesize\(u_x\)}
   \uput{3}[180](10,60){\footnotesize\(u_y\)}
   \uput{3}[0](210,60){\footnotesize\(u_{y'}\)}

   \rput(10,20){\secunder{R_{xy}}{\la_{xy}}}
   \rput(110,20){\secunderf{R_{xy'}}{\la_{xy'}}}
   \uput{3}[270](110,16){\footnotesize\(u_x\)}
   \uput{3}[180](10,20){\footnotesize\(u_y\)}
   \uput{3}[0](210,20){\footnotesize\(u_{y'}\)}
  \end{pspicture}
 \end{center}
 where the second line shows an equivalent rewriting of the original setting. However, with this secunder representation our formula follows easily after realising that
 \begin{itemize}
  \item the potentials on the \(x\)-side of the amplifiers are \(\la_{xy}\cdot u_y\) for the respective vertices \(y\);
  \item from here the potential \(u_x\) is computed using the well-known formula for a voltage divider (of conductances \(C^\text{se}_{xy}=1/R^\text{se}_{xy}\)).
 \end{itemize}
 Putting all that in formulas, we have
 \[
  \ba
   u_x&=\sum_{\substack{y\sim x\\y\ne x}}\la_{xy}u_y\cdot\frac{C^\text{se}_{xy}}{\sum\limits_{\substack{z\sim x\\z\ne x}}C^\text{se}_{xz}}=\sum_{\substack{y\sim x\\y\ne x}}u_y\frac{C_{xy}\frac{2\la_{xy}}{\la_{xy}+1}}{\sum\limits_{\substack{z\sim x\\z\ne x}}C_{xz}\frac2{\la_{xz}+1}}=\frac{\sum\limits_{y\sim x}u_yD_{xy}\ga_{xy}-u_xD_{xx}}{D_x-D_{xx}}
  \ea
 \]
 with the use of \eqref{eq:prisech}, \eqref{eq:ddef} and \(\la_{xx}=1\). Rearranging finishes the proof.
\end{proof}

Let \(P_{\cdot\cdot}\) be the transition probabilities of an irreducible Markov chain on the finite, connected graph \(G\). Throughout this manuscript we assume that \(P_{xy}>0\) whenever \((x,\,y)\in E\) is an edge of the graph, totally asymmetric steps are not handled by our methods (although we suspect that a meaningful limit could be worked out for these cases). As usual, the graph has a loop on vertex \(x\) whenever \(P_{xx}>0\). We now give a recipe of how to build an electric network for this chain so that the resulting voltages and currents have the classical probabilistic interpretations (see e.g., Doyle-Snell \cite{doyle_snell}). This will be done regardless of whether the chain is reversible or irreversible. The unique stationary distribution of the chain will be called \(\mu_\cdot\), and we now make the following choices:
\be
 \ba
  D_{xy}:&=\sqrt{\mu_x\cdot P_{xy}\cdot\mu_y\cdot P_{yx}};\\
  \ga_{xy}:&=\sqrt{\frac{\mu_x\cdot P_{xy}}{\mu_y\cdot P_{yx}}}.
 \ea\label{eq:ch2circ}
\ee
Notice first that these choices are consistent with the respective symmetry and log-antisymmetry of \(D_{\cdot\cdot}\) and \(\ga_{\cdot\cdot}\). It is also clear that the conductances \(C_{\cdot\cdot}\) and the amplifying factors \(\la_{\cdot\cdot}\) can also be expressed with the help of the above quantities. Following the definition \eqref{eq:ddef}, we have
\be
 D_x=\sum_{z\sim x}D_{xz}\ga_{zx}=\sum_{z\sim x}\mu_zP_{zx}=\mu_x.\label{eq:dmu}
\ee
Recall the non-intersecting subsets \(A\) and \(B\) of the vertex set \(V\), and define, for \(x\in V\), the first reaching times
\be
 \tau^0_A:\,=\inf\{t\ge0\,:\,X(t)\in A\}\label{eq:tau0}
\ee
and similarly \(\tau^0_B\) of these sets by the Markov chain started from \(x\). When \(x\in A\) (\(B\)), we define \(\tau^0_A\) (\(\tau^0_B\), respectively) to be 0. \(\Pv_x\) will stand for the probabilities associated with the chain started from \(x\). For short, we will set \(h_x:\,=\Pv_x\{\tau^0_A<\tau^0_B\}\).
\begin{tm}\label{tm:base}
 Set up the electric network with the choices \eqref{eq:ch2circ}, apply constant potentials \(U_A\equiv1\) on vertices of the set \(A\), \(U_B\equiv0\) on vertices of the set \(B\), and make no external connections to vertices of \(V-A-B\). Then for every \(x\in V\) we have \(u_x=h_x\).
\end{tm}
\begin{proof}
 By definition we have \(h_x=1\) for vertices \(x\in A\), \(h_x=0\) for vertices \(x\in B\), and by a first step analysis of the Markov chain,
 \[
  h_x=\sum_{y\sim x}h_yP_{xy}
 \]
 when \(x\in V-A-B\). Next we find that by definition we have \(u_x=1\) for vertices \(x\in A\), \(u_x=0\) for vertices \(x\in B\), and by \eqref{eq:uuu}, \eqref{eq:ch2circ} and \eqref{eq:dmu},
 \[
  u_x=\sum_{y\sim x}u_y\frac{D_{xy}\gamma_{xy}}{D_x}=\sum_{y\sim x}u_yP_{xy}
 \]
 when \(x\in V-A-B\). Thus \(h_\cdot\) and \(u_\cdot\) satisfy the same (well defined) equations with the same boundary conditions, therefore they agree on all vertices.
\end{proof}

A nice consequence of the analogy is what happens to our electric network when we reverse our Markov chain. The reversed Markov chain has the same stationary distribution \(\mu_x\) as the original one, and its transition probabilities become
\[
 \hat P_{xy}=P_{yx}\cdot\frac{\mu_y}{\mu_x}.
\]
We will simply call the network that corresponds the reversed chain the \emph{reversed network}, and its parameters will be marked by hats. They are
\be
 \ba
  \hat D_{xy}&=\sqrt{\mu_x\cdot\hat P_{xy}\cdot\mu_y\cdot\hat P_{yx}}=\sqrt{\mu_y\cdot P_{yx}\cdot\mu_x\cdot P_{xy}}=D_{xy};\\
  \hat\ga_{xy}&=\sqrt{\frac{\mu_x\cdot\hat P_{xy}}{\mu_y\cdot\hat P_{yx}}}=\sqrt{\frac{\mu_y\cdot P_{yx}}{\mu_x\cdot P_{xy}}}=\ga_{yx}=\frac1{\ga_{xy}}.
 \ea\label{eq:revddef}
\ee
This also implies \(\hat C_{xy}=C_{xy}\) and \(\hat\la_{xy}=\la_{yx}=1/\la_{xy}\) for the reversed network, in other words \emph{reversing the Markov chain simply reverses the direction of our voltage amplifiers while keeps the resistance values intact}. A Markov chain is reversible if and only if the corresponding network has all its amplifiers with \(\la_{xy}\equiv1\). Indeed an amplifier of parameter 1 is just a plain wire, therefore this case reduces to the classical reversible setting with ordinary resistors on the edges.

\subsection{Markovian networks}

Notice that \(P_{\cdot\cdot}\) being a Markov transition probability imposes restrictions on our electric network. From now on, \(\sum\limits_{z\sim x\in V}\) will be our notation for double summation on all neighbouring vertices \(x\) and \(z\) in \(V\).
\begin{tm}
 Suppose that we are given an electric network built of our components on the edges of the finite connected graph \(G=(V,\,E)\). There is an irreducible Markov chain of graph \(G\) and transition probabilities \(P_{\cdot\cdot}\) such that \eqref{eq:ch2circ} holds if and only if we have both
 \begin{align}
  \sum_{z\sim x}D_{xz}\ga_{xz}&=\sum_{z\sim x}D_{xz}\ga_{zx}\qquad(\forall x\in V)\text{, and}\label{eq:dx}\\
  \sum_{z\sim x\in V}D_{xz}\ga_{zx}&=1.\label{eq:dnorm}
 \end{align}
 In this case
 \be
  \mu_x=\sum_{z\sim x}D_{xz}\ga_{zx}\,\qquad\text{and}\qquad P_{xy}=\frac{D_{xy}\ga_{xy}}{\sum\limits_{z\sim x}D_{xz}\ga_{zx}}.\label{eq:circ2ch}
 \ee
\end{tm}
\begin{proof}
 If \eqref{eq:ch2circ} holds for a Markov transition probability \(P_{\cdot\cdot}\), then the above formulas follow from direct verification. Conversely, if \eqref{eq:dx} and \eqref{eq:dnorm} hold then we make the definition
 \[
  P_{xy}=\frac{D_{xy}\ga_{xy}}{\sum\limits_{z\sim x}D_{xz}\ga_{zx}},
 \]
 and notice that
 \[
  \sum_{y\sim x}P_{xy}=\sum_{y\sim x}\frac{D_{xy}\ga_{xy}}{\sum\limits_{z\sim x}D_{xz}\ga_{zx}}=1\qquad(\forall x\in V),
 \]
 which shows that \(P_{\cdot\cdot}\) is a Markov transition probability matrix. It is irreducible by positivity of our parameters, and its stationary distribution \(\mu_\cdot\) is the unique vector with
 \[
  \ba
   \mu_y&=\sum_{x\sim y}\mu_xP_{xy}=\sum_{x\sim y}\mu_x\frac{D_{xy}\ga_{xy}}{\sum\limits_{z\sim x}D_{xz}\ga_{zx}}\qquad(\forall y\in V)\text{ and}\\
   \sum_{x\in V}\mu_x&=1.
  \ea
 \]
 Notice, however, that the vector \(\bigl(\sum\limits_{z\sim x}D_{xz}\ga_{zx}\bigr)_{x\in V}\) satisfies the same properties:
 \[
  \sum_{x\sim y}\sum_{z\sim x}D_{xz}\ga_{zx}P_{xy}=\sum_{x\sim y}\sum_{z\sim x}D_{xz}\ga_{zx}\frac{D_{xy}\ga_{xy}}{\sum\limits_{w\sim x}D_{xw}\ga_{wx}}=\sum_{x\sim y}D_{xy}\ga_{xy}=\sum_{x\sim y}D_{yx}\ga_{xy}
 \]
 by the symmetry of \(D_{\cdot\cdot}\), and
 \[
  \sum_{z\sim x\in V}D_{xz}\ga_{zx}=1
 \]
 by \eqref{eq:dnorm}. Therefore these two vectors agree.
\end{proof}
\begin{rem}
 The normalisation \eqref{eq:dnorm} is just an artificial choice, and is not essential at all. Given an electric network, multiplying every resistor value by the same constant \(K\) while keeping the amplifiers unchanged will result in the same voltages everywhere with currents multiplied by \(1/K\). In particular, Theorem \ref{tm:base} holds true in this case.
\end{rem}
\begin{rem}
 The condition \eqref{eq:dx} is, on the other hand, very essential, and we will refer to networks with this property as \emph{Markovian}. On a technical level it states that we can extend the definition \eqref{eq:dmu} of \(D_\cdot\) by
 \[
  D_x=\sum_{z\sim x}D_{xz}\ga_{zx}=\sum_{z\sim x}D_{xz}\ga_{xz}.
 \]
 Notice also that this implies
 \be
  D_x=\sum_{z\sim x}D_{xz}\frac{\ga_{zx}+\ga_{xz}}2=\sum_{z\sim x}C_{xz}.\label{eq:cd}
 \ee
 The Markovian property also has a rather intuitive meaning: considering \eqref{eq:uuu} it states that the constant potential \(u_x\equiv U\) for all vertices is a valid solution of the (free) network.
\end{rem}
This was, of course, trivially true for the all-resistors networks that correspond to reversible Markov chains. Consider a set of connected resistors, and apply potential \(U\) on one of the vertices. Then all vertices will stay at potential \(U\), with no current flowing anywhere in the network. This is not at all straightforward with our generalised networks of resistors and amplifiers. Applying potential \(U\) on one of the vertices, the amplifiers will change voltages for different parts of the network, and this can keep up currents in the cycles of the graph \(G\). The Markovian property is that, nevertheless, each vertex will still stay at the same potential \(U\) even if circular (that is, divergence free) currents flow in the system.

A classical result for Markov chains follows easily from the analogy.
\begin{cor}
 A Markov chain is reversible if and only if for every closed cycle \(x_0,\,x_1,\,x_2,\,\dots,\,x_n=x_0\) in the graph \(G\) we have
 \[
  P_{x_0x_1}\cdot P_{x_1x_2}\cdots P_{x_{n-1}x_0}=P_{x_0x_{n-1}}\cdot P_{x_{n-1}x_{n-2}}\cdots P_{x_1x_0}.
 \]
 In particular, any Markov chain on a finite connected tree \(G\) is necessarily reversible.
\end{cor}
\begin{proof}
 Rewriting the above formula and using \eqref{eq:circ2ch} together with the symmetry of \(D_{\cdot\cdot}\), we arrive to the equivalent statement
 \[
  \ba
   &\ga_{x_0x_1}\cdot\ga_{x_1x_2}\cdots\ga_{x_{n-1}x_0}=\ga_{x_0x_{n-1}}\cdot\ga_{x_{n-1}x_{n-2}}\cdots\ga_{x_1x_0}\text{, or}\\
   &\la_{x_0x_1}\cdot\la_{x_1x_2}\cdots\la_{x_{n-1}x_0}=1.
  \ea
 \]
 This is of course trivially true in the reversible case where all of the amplifiers have \(\la_{xy}=1\). For the other direction, assume now the above formula to hold, and turn it into electrical language. It says that the total multiplication factor of the potentials is one along any closed cycle of the circuit. It follows that fixing one vertex at potential \(U\), zero currents everywhere in the network is a solution. By uniqueness this is the only solution. The network being Markovian on the other hand tells us that every vertex has to be on potential \(U\). With no currents the only way this can happen is that all of the amplifiers have parameter one, and the chain is reversible.

 A similar argument works directly for the tree: since there are no cycles, no current can flow if only one vertex is fixed at potential \(U\). The Markovian property again tells us that every vertex will be on potential \(U\) which again means parameter one for all of the amplifiers and thus reversibility of the chain.
\end{proof}

\subsection{Effective resistance, capacity, and escape probabilities}

In this section we make sense of effective resistance in our network, and give it a probabilistic interpretation similar to that of the classical case. The setting is the one of Proposition \ref{pr:uuu}: two disjoint subsets \(A\) and \(B\) of the vertices are forced to be on constant potentials \(U_A\) and \(U_B\), respectively. Define the \emph{total incoming current to the set \(A\)} (c.f.\ \eqref{eq:ic}) as
\be
i_A:\,=\sum_{x\in A}i_x=\sum_{y\sim x\in A}i_{xy}.\label{eq:tic}
\ee
Notice that by conservation of currents, this agrees to the sum of currents of edges across the boundary of \(A\), and also \(i_A+i_B=0\). The existence of the effective resistance between sets \(A\) and \(B\) means that the network between these sets can be replaced by a single resistor. This is not true for arbitrary configurations, since the amplifiers in general push the characteristics away from that of a single resistor. It is, however, true for networks that match a Markov chain, this is formulated in the next theorem:
\begin{tm}
 In a Markovian electric network, for any disjoint \(A,\,B\subset V\) there is a constant \(\reff_{AB}>0\) such that
 \[
  U_A-U_B=\reff_{AB}\cdot i_A\qquad(\forall U_A,\,U_B\in\Rb).
 \]
\end{tm}
\begin{proof}
 The proof will again proceed along the lines of linearity. When \(U_A=U_B\) then we just have the Markovian solution with zero incoming currents, thus \(i_A=0\) and everything is trivial. Suppose that we are given arbitrary reals \(U_A\ne U_B,\,U'_A\ne U'_B\).
 \begin{itemize}
  \item We consider two solutions of our network: the one \(u_\cdot\), \(i_{\cdot\cdot}\) that satisfies the given boundary conditions \(u_x\equiv U_A\) for \(x\in A\) and \(u_x\equiv U_B\) for \(x\in B\), and one that comes from the Markovian property: \(u^\text M_x\equiv U_B\), \(i^\text M_x\equiv0\) (incoming currents to vertex \(x\), not to be mixed with currents \(i^\text M_{xy}\) of edges!) for all \(x\in V\). We think of this latter one as a solution with boundary conditions \(u^\text M_x\equiv U_B\) for all \(x\in A\cup B\).
  \item The difference \(u-u^\text M\) and \(i-i^\text M\) of these two is yet another solution due to linearity. It has boundary conditions \(u_x-u^M_x\equiv U_A-U_B\) on \(x\in A\), \(u_x-u^M_x\equiv U_B-U_B=0\) on \(x\in B\), and notice that the incoming current to the set \(A\) is still \(i_A-i^\text M_A=i_A-0=i_A\).
  \item Again by linearity every current and potential can be multiplied by the factor \(\frac{U'_A-U'_B}{U_A-U_B}\) and we still have a valid system. This looks like
  \[
   (u_\cdot-u^\text M_\cdot)\cdot\frac{U'_A-U'_B}{U_A-U_B},\qquad(i_{\cdot\cdot}-i^\text M_{\cdot\cdot})\cdot\frac{U'_A-U'_B}{U_A-U_B}
  \]
  and therefore has boundary conditions
  \[
   \ba
    &(U_A-U_B)\cdot\frac{U'_A-U'_B}{U_A-U_B}=U'_A-U'_B\text{ on }x\in A,\qquad\text{and}\\
    &0\cdot\frac{U'_A-U'_B}{U_A-U_B}=0\text{ on }x\in B,
   \ea
  \]
  with incoming current \(i_A\cdot\frac{U'_A-U'_B}{U_A-U_B}\) on the set \(A\).
  \item Finally, add the Markovian solution with constant potential \({u_\cdot^\text M}'\equiv U'_B\) everywhere, and \({i^\text M_{\cdot\cdot}}'\) with zero incoming currents in all vertices. This results in
  \[
   (u_\cdot-u_\cdot^\text M)\cdot\frac{U'_A-U'_B}{U_A-U_B}+U'_B,\qquad(i_{\cdot\cdot}-i_{\cdot\cdot}^\text M)\cdot\frac{U'_A-U'_B}{U_A-U_B}+{i^\text M}'
  \]
  and therefore has boundary conditions
  \[
   \ba
    &U'_A-U'_B+U'_B=U'_A\text{ on }x\in A,\qquad\text{and}\\
    &0+U'_B=U'_B\text{ on }x\in B,
   \ea
  \]
  with incoming current \(i_A\cdot\frac{U'_A-U'_B}{U_A-U_B}+0=i_A\cdot\frac{U'_A-U'_B}{U_A-U_B}\) on the set \(A\).
 \end{itemize}
 We have thus produced the solution for the boundary conditions \(U_A'\) and \(U_B'\), and concluded that the corresponding incoming current to the set \(A\) is
 \[
  i_A'=i_A\cdot\frac{U'_A-U'_B}{U_A-U_B}.
 \]
 This is equivalent to the statement of the theorem, i.e., the ratio \((U_A-U_B)/i_A\) is a constant for all boundary potentials \(U_A\) and \(U_B\).
\end{proof}

A Markovian network has a further peculiar property: the effective resistance \(\reff_{AB}\) stays the same if we reverse each of the amplifiers. Recall that the network then turns into that of the reversed Markov chain. To prove this property, we follow Slowik's argument \cite{slowik_var}, and first define the one-step Markov generator \(L\) on functions \(f\,:\,V\to\Rb\):
\[
 (Lf)_x=\sum_{y\sim x}P_{xy}(f_y-f_x).
\]
Rewriting this via \eqref{eq:circ2ch} into electrical terms and then applying \eqref{eq:dx} we get
\[
 -(Lf)_x=\sum_{y\sim x}\frac{D_{xy}\ga_{xy}}{D_x}(f_x-f_y)=\frac1{D_x}\sum_{y\sim x}D_{xy}\cdot(\ga_{yx}f_x-\ga_{xy}f_y).
\]
This latter formula is meaningful in electrical terms, as soon as we imagine \(f_\cdot\) as a potential applied on vertices of the graph, and define the resulting currents
\be
i^f_{xy}:\,=D_{xy}\cdot(\ga_{yx}f_x-\ga_{xy}f_y)\label{eq:iu}
\ee
via \eqref{eq:idg}. We thus see, c.f.\ \eqref{eq:ic}, that
\[
 -(Lf)_x=\frac{i^f_x}{D_x}.
\]
\(i^f_x\) is the current we are required to pump in vertex \(x\) in order to maintain potential \(f_x\). The quantity
\[
 \mathcal E(f):\,=\sum_{x\in V}\mu_xf_x\cdot(-Lf)_x=\sum_xf_x\cdot i^f_x
\]
is referred to as the \emph{energy} associated to the pair \(P_{\cdot\cdot},\,\mu_\cdot\), and we now see that it is the total electric power we need to pump in the system in order to maintain potential \(f_x\) at each vertex \(x\). (As usual, we do not count the external energy sources (absorbers) required by the amplifiers to work.)

With this preparation we now prove
\begin{pr}\label{pr:revres}
 Reversing a Markovian network does not affect the effective resistance:
 \[
  \widehat{\reff_{AB}}=\reff_{AB}.
 \]
\end{pr}
\begin{proof}
 We repeat Slowik's arguments \cite{slowik_var} in the electrical language. Take two functions \(f\) and \(g\) on \(V\), and apply \eqref{eq:dx} in the first term and symmetry of the double summation and of \(D_{\cdot\cdot}\) in the second term below:
 \be
  \ba
   \sum_xf_x\cdot i^g_x&=\sum_{x\in V}f_x\sum_{y\sim x}D_{xy}\cdot(\ga_{yx}g_x-\ga_{xy}g_y)\\
   &=\sum_{x\in V}f_xg_x\sum_{y\sim x}D_{xy}\ga_{yx}-\sum_{y\sim x\in V}f_xD_{xy}\ga_{xy}g_y\\
   &=\sum_{x\in V}f_xg_x\sum_{y\sim x}D_{xy}\ga_{xy}-\sum_{y\sim x\in V}g_xD_{xy}\ga_{yx}f_y\\
   &=\sum_{x\in V}g_x\sum_{y\sim x}D_{xy}\cdot(\hat\ga_{yx}f_x-\hat\ga_{xy}f_y)=\sum_xg_x\cdot\hat i^f_x.
  \ea\label{eq:adj}
 \ee
 (This equation is the electrical way of saying that the adjoint of the generator is the one of the reversed process.) As before, fix the boundary conditions \(u_x\equiv1\equiv\hat u_x\) on \(x\in A\) and \(u_x\equiv0\equiv\hat u_x\) on \(x\in B\) for two scenarios: \(u_\cdot,\,i_{\cdot\cdot}\) of the original network and \(\hat u_\cdot,\,\hat i_{\cdot\cdot}\) of the reversed one. This latter has all its amplifiers reversed, and it corresponds to the reversed Markov chain. We claim that in our situation \(i^u_x\equiv0\equiv\hat i^{\hat u}_x\) on \(x\notin A\cup B\), since these are free vertices. This, together with the common boundary condition for the two networks implies
 \[
  \ba
   \mathcal E(u)&=\sum_{x\in V}u_x\cdot i^u_x=\sum_{x\in A}u_x\cdot i^u_x=\sum_{x\in A}\hat u_x\cdot i^u_x=\sum_{x\in V}\hat u_x\cdot i^u_x\\
   &=\sum_{x\in V}u_x\cdot\hat i^{\hat u}_x=\sum_{x\in A}u_x\cdot\hat i^{\hat u}_x=\sum_{x\in A}\hat u_x\cdot\hat i^{\hat u}_x=\sum_{x\in V}\hat u_x\cdot\hat i^{\hat u}_x=\hat{\mathcal E}(\hat u).
  \ea
 \]
 Rewriting the power we apply to maintain our boundary conditions gives
 \[
  \ceff_{AB}=(U_A-U_B)^2\cdot\ceff_{AB}=\mathcal E(u)=\hat{\mathcal E}(\hat u)=(U_A-U_B)^2\cdot\widehat{\ceff_{AB}}=\widehat{\ceff_{AB}}.
 \]
\end{proof}

Next we introduce what is called the \emph{capacity} in the theory of Markov chains, and show that it has close connections to the effective resistance. We again assume that \(A,\,B\subset V\) are non-empty, disjoint, and follow Gaudilli\`ere-Landim and Slowik \cite{gau-lan,slowik_var} by defining
\[
\tau_A:\,=\inf\{t>0\,:\,X(t)\in A\}
\]
(c.f.\ \eqref{eq:tau0}) and
\[
 \cp(A,\,B):\,=\sum_{x\in A}\mu_x\Pv_x\{\tau_B<\tau_A\}.
\]
\begin{pr}
 The above capacity is simply the effective conductance \(\ceff_{AB}=1/\reff_{AB}\) between the sets \(A\) and \(B\).
\end{pr}
\begin{proof}
 We use the analogy set up in Theorem \ref{tm:base} as
 \[
  \ba
   \cp(A,\,B)&=\sum_{x\in A}\mu_x\sum_{y\sim x}P_{xy}\cdot\Pv_y\{\tau^0_B<\tau^0_A\}\\
   &=\sum_{y\sim x\in A}\mu_xP_{xy}\cdot(1-\Pv_y\{\tau^0_A<\tau^0_B\})\\
   &=\sum_{y\sim x\in A}D_{xy}\ga_{xy}(1-u_y)=\sum_{y\sim x\in A}D_{xy}(\ga_{yx}\cdot1-\ga_{xy}u_y)\\
   &=\sum_{y\sim x\in A}D_{xy}(\ga_{yx}u_x-\ga_{xy}u_y)=\sum_{y\sim x\in A}i_{xy}\\
   &=i_A=U_A\cdot\ceff_{AB}=\ceff_{AB}.
  \ea
 \]
Along the way we also used \eqref{eq:dx}, the fact that \(u_x\equiv U_A=1\) for all \(x\in A\), and finally \eqref{eq:idg}.
\end{proof}
It follows immediately that the capacity is a symmetric quantity in its two arguments \(A\) and \(B\). The identity \(\cp(A,\,B)=\widehat\cp(B,\,A)\) also follows from the previous proposition.
\begin{rem}
 Gaudilli\`ere-Landim and Slowik \cite{gau-lan,slowik_var} also establish
 \[
  \cp(A,\,B)=\frac12\sum_{x\sim y\in V}\mu_xP_{xy}^\text s(h_x-h_y)^2,
 \]
 with the symmetrised transitions \(P_{xy}^\text s=\frac12(P_{xy}+\hat P_{xy})\). \eqref{eq:circ2ch} together with \eqref{eq:revddef} and \eqref{eq:ddef} gives
 \be
  P_{xy}^\text s=\frac{D_{xy}\ga_{xy}+D_{xy}\ga_{yx}}{2D_x}=\frac{C_{xy}}{D_x},\label{eq:ps}
 \ee
 and the capacity gets another interesting interpretation: in the setting of Theorem \ref{tm:base},
 \be
  \cp(A,\,B)=\frac12\sum_{x\sim y\in V}C_{xy}(u_x-u_y)^2,\label{eq:cape}
 \ee
 \emph{the ohmic power loss on the resistors, should we apply the actual voltages \(u_\cdot\) on them without the amplifiers}. It is important to note that this interpretation is non-physical: with the amplifiers the ohmic losses are not given by the above formula, without the amplifiers the voltages \(u_\cdot\) would be totally different.
\end{rem}
The capacity, being \(\ceff_{AB}\) is, however, equal to the total power \((U_A-U_B)^2\cdot\ceff_{AB}=(1-0)^2\cdot\ceff_{AB}\) we need to pump in the set \(A\) to keep it on potential \(U_A=1\).

We repeat the computation for \eqref{eq:cape} in the electrical language. First notice that by \eqref{eq:iu} and \eqref{eq:ddef},
\[
 \frac12(i^u_{xy}+\hat i^u_{xy})=D_{xy}\cdot\Bigl(\frac{\ga_{yx}+\ga_{xy}}2u_x-\frac{\ga_{xy}+\ga_{yx}}2u_y\Bigr)=C_{xy}(u_x-u_y),
\]
fixing the voltages everywhere (!), the average of the current and the reversed current is the one of the network without amplifiers. This is the starting point to expand the right hand-side of \eqref{eq:cape}:
\begin{multline*}
\frac12\sum_{x\sim y\in V}C_{xy}(u_x-u_y)^2\\
 \ba
  &=\frac14\sum_{x\sim y\in V}(u_x-u_y)(i^u_{xy}+\hat i^u_{xy})\\
  &=\frac14\sum_{x\sim y\in V}u_xi^u_{xy}+\frac14\sum_{x\sim y\in V}u_x\hat i^u_{xy}-\frac14\sum_{x\sim y\in V}u_yi^u_{xy}-\frac14\sum_{x\sim y\in V}u_y\hat i^u_{xy}\\
  &=\frac14\sum_{x\sim y\in V}u_xi^u_{xy}+\frac14\sum_{x\sim y\in V}u_x\hat i^u_{xy}+\frac14\sum_{x\sim y\in V}u_yi^u_{yx}+\frac14\sum_{x\sim y\in V}u_y\hat i^u_{yx}\\
  &=\frac14\sum_{x\in V}u_xi^u_x+\frac14\sum_{x\in V}u_x\hat i^u_x+\frac14\sum_{y\in V}u_yi^u_y+\frac14\sum_{y\in V}u_y\hat i^u_y\\
  &=\frac14\sum_{x\in V}u_xi^u_x+\frac14\sum_{x\in V}u_xi^u_x+\frac14\sum_{y\in V}u_yi^u_y+\frac14\sum_{y\in V}u_yi^u_y=\sum_{x\in V}u_xi^u_x.
 \ea
\end{multline*}
with the use of the adjoint identity \eqref{eq:adj}. Now, as in the proof of Proposition \ref{pr:revres}, apply our usual boundary conditions, and the right hand-side becomes the total power required to maintain the boundary conditions or, equivalently, \(\ceff_{AB}\).

Finally, we define the escape probability from the set \(A\) and show its connection to the effective resistance, this goes exactly as in the reversible case. Suppose that the Markov chain is started from its stationary distribution \(\mu\), conditioned on being in the set \(A\). (When \(A={a}\) is a singleton, this is just the unit mass on the vertex \(a\).) The escape probability is the chance that the chain reaches set \(B\) before its first return to \(A\):
\[
 \Pv\{\tau_B<\tau_A\}=\sum_{x\in A}\frac{\mu_x}{\mu(A)}\Pv_x\{\tau_B<\tau_A\}=\frac{\cp(A,\,B)}{\mu(A)}=\frac{\ceff_{AB}}{\sum\limits_{z\in A}D_z}=\frac{\ceff_{AB}}{\sum\limits_{y\sim z\in A}C_{zy}}.
\]
The last step used \eqref{eq:cd}. The right hand-side agrees word for word with the classical reversible result, the starting point of elegant recurrence-transience proofs.

By symmetrising a Markov chain we mean replacing its transition probabilities by \(P^\text s_{\cdot\cdot}\).
\begin{cor}
 Symmetrising a Markov chain never increases the escape probabilities.
\end{cor}
\begin{proof}
 By symmetrising we keep the stationary distribution \(\mu_\cdot=D_\cdot\) and the conductances \(C_{\cdot\cdot}\) unchanged while the amplifiers all become trivial: \(\la_{\cdot\cdot}\equiv1\). This can easily be seen via \eqref{eq:revddef} and \eqref{eq:ps}. Denote the potentials that result our usual boundary conditions \(U_A\equiv1\) and \(U_B\equiv0\) by \(u_\cdot\) in the original network and by \(u^\text s_\cdot\) in the symmetrised one. The classical Dirichlet principle for the reversible case tells us that \(u^\text s_\cdot\) is the potential that minimises the ohmic power losses in the resistors for the reversible chain. Therefore
 \[
  \frac12\sum_{x\sim y\in V}C_{xy}(u^\text s_x-u^\text s_y)^2\le\frac12\sum_{x\sim y\in V}C_{xy}(u_x-u_y)^2.
 \]
 Since the conductances agree for the two networks, the left hand-side is the capacity of the symmetrised network, while the right hand-side is the one of the original network. Dominance of the capacities then implies that of the escape probabilities as the conductances are not changed by symmetrising.
\end{proof}

\subsection{The current}

In this section we give a probabilistic interpretation of the currents \(i_{\cdot\cdot}\) of the network. We take a singleton set \(A=\{a\}\) and another, arbitrary set \(B\not\ni a\). Start the Markov chain from \(a\), and define \(v_x\) as the expected number of visits to vertex \(x\) before the first hitting of the set \(B\). Clearly \(v_x\equiv0\) for all \(x\in B\). A last step analysis shows that for any \(a\ne x\notin B\) we have
\[
 v_x=\sum_{y\sim x}v_yP_{yx}=\sum_{y\sim x}v_y\hat P_{xy}\cdot\frac{\mu_x}{\mu_y},
\]
in other words \(v_x/\mu_x\) is harmonic w.r.t.\ \(\hat P\). The boundary conditions are \(v_a/\mu_a\) is fixed (to be determined later), and \(v_x/\mu_x\equiv0\) on \(x\in B\). Consider now the corresponding electric network of parameters \(C_{\cdot\cdot}\) and \(\hat\la_{\cdot\cdot}\) with these boundary conditions. It follows that its potentials \(\hat u_x\) agree with \(v_x/\mu_x\) for all \(x\in V\), and its currents are (see \eqref{eq:idg})
\[
 \hat i_{xy}=D_{xy}\cdot\bigl(\hat\ga_{yx}\cdot\frac{v_x}{\mu_x}-\hat\ga_{xy}\cdot\frac{v_y}{\mu_y}\bigr)=D_{xy}\cdot\bigl(\ga_{xy}\cdot\frac{v_x}{\mu_x}-\ga_{yx}\cdot\frac{v_y}{\mu_y}\bigr)=v_xP_{xy}-v_yP_{xy}.
\]
This latter is the expected number of jumps from \(x\) to \(y\) minus that from \(y\) to \(x\) before absorption in \(B\). It remains to fix the boundary term \(\hat u_a=v_a/\mu_a\). This is done by the simple observation that the chain has to exit vertex \(a\) one more times than enter it, thus
\[
 \hat i_a=\sum_{y\sim a}\hat i_{ay}=1.
\]
This fixes a unique potential \(\hat u_a=v_a/\mu_a\) for boundary. Everything is analogue to the classical reversible case, except that we need to use the reversed network.

\subsection{Commute times and costs}

We have found that a classical result of Chandra, Raghavan, Ruzzo, Smolensky and Tiwari \cite{commt-res} on effective resistance and commute times (or costs) can be extended easily to the irreversible case. Work on commute times in this case has also been done by Doyle and Steiner \cite{doyle_steiner_comm_geom}. Fix two vertices \(a\ne b\) of the graph and a cost function \(k_{xy}\) on edges of the graph. Costs \(k_{xy}\) and \(k_{yx}\) can be different, we do not require any relation between these two. The expected cost of the chain from \(a\) to \(b\) is the expected price to pay until the first hitting of the chain to vertex \(b\) if started from \(a\). It is defined as
\[
 H^k_{ab}:\,=\Ev_a\sum_{t=1}^{\tau_b^0}k_{X(t-1)\,X(t)},
\]
where we (ab)used definition \eqref{eq:tau0} (by writing \(\tau^0_b\) for \(\tau^0_{\{b\}}\)). We consider an empty sum to be zero for the case \(H^k_{aa}=0\). In particular, for \(k\equiv1\) we arrive to the expected hitting time \(H^1_{ab}=\Ev_a(\tau^0_b)\). Define the expected \emph{commute cost} \(K^k_{ab}=H^k_{ab}+H^k_{ba}\), this becomes the expected \emph{commute time} for \(k\equiv1\).
\begin{tm}
 The expected commute cost can be computed by
 \[
  K_{ab}^k=\reff_{ab}\cdot D^k,
 \]
  with \(D^k:\,=\sum\limits_{y\sim x\in V}D_{xy}\ga_{xy}k_{xy}\).
\end{tm}
We spell out the case \(k\equiv1\): the expected commute time is
\[
 K_{ab}^1=\reff_{ab}\cdot D^1=\reff_{ab}\cdot\sum_{y\sim x\in V}D_{xy}\ga_{xy}=\reff_{ab}\cdot\sum_{x\in V}D_x=\reff_{ab}\cdot\sum_{z\sim x\in V}C_{xz}
\]
via \eqref{eq:cd}. This is the exact same formula as the one of \cite{commt-res} for the reversible case.
\begin{proof}
 We start with a first step analysis and write, for any \(x\ne b\),
 \be
 H_{xb}^k=\sum_{y\sim x}P_{xy}(k_{xy}+H_{yb}^k)=\sum_{y\sim x}\frac{D_{xy}\ga_{xy}}{D_x}(k_{xy}+H_{yb}^k).\label{eq:heq}
 \ee
 In our electric network, impose boundary conditions \(u_x=H^k_{xb}\) on each vertex \(x\). This results in currents
 \[
  i_{xy}=D_{xy}\cdot(\ga_{yx}H_{xb}^k-\ga_{xy}H^k_{yb})
 \]
 according to \eqref{eq:idg}, and, via \eqref{eq:heq}, the necessity of pumping external currents
 \[
  \ba
   i_x=\sum_{y\sim x}i_{xy}&=\sum_{y\sim x}D_{xy}\ga_{yx}H_{xb}^k-\sum_{y\sim x}D_{xy}\ga_{xy}H^k_{yb}\\
   &=D_xH_{xb}^k-D_xH_{xb}^k+\sum_{y\sim x}D_{xy}\ga_{xy}k_{xy}\\
   &=\sum_{y\sim x}D_{xy}\ga_{xy}k_{xy}
  \ea
 \]
 into each vertex \(x\ne b\). By conservation of current,
 \[
  i_b=\sum_{y\sim b}D_{by}\ga_{by}k_{by}-D^k
 \]
 with \(D^k:\,=\sum\limits_{y\sim x\in V}D_{xy}\ga_{xy}k_{xy}\).

 A second configuration we consider is \(u'_x=H^k_{xa}\) on each vertex \(x\), in a similar fashion this has external currents
 \[
  i'_x=\sum_{y\sim x}D_{xy}\ga_{xy}k_{xy}
 \]
 for all \(x\ne a\), and
 \[
  i'_a=\sum_{y\sim a}D_{ay}\ga_{ay}k_{ay}-D^k.
 \]
 Our equations being linear, the difference \(u-u'\) is also a solution of the network. It has potentials and external currents
 \[
  \ba
   u_a-u'_a&=H^k_{ab}-H^k_{aa}=H^k_{ab}&&\text{and}\\
   i_a-i'_a&=\sum_{y\sim a}D_{ay}\ga_{ay}k_{ay}-\sum_{y\sim a}D_{ay}\ga_{ay}k_{ay}+D^k=D^k\quad&&\text{in vertex }a,\\
   u_b-u'_b&=H^k_{bb}-H^k_{ba}=-H^k_{ba}&&\text{and}\\
   i_b-i'_b&=\sum_{y\sim b}D_{by}\ga_{by}k_{by}-D^k-\sum_{y\sim b}D_{by}\ga_{by}k_{by}=-D^k\quad&&\text{in vertex }b,\\
   i_x-i'_x&=\sum_{y\sim x}D_{xy}\ga_{xy}k_{xy}-\sum_{y\sim x}D_{xy}\ga_{xy}k_{xy}=0\quad&&\text{elsewhere}.
  \ea
 \]
 Therefore this combination only has boundary conditions at \(a\) and \(b\), all other vertices are free. The effective conductance between \(a\) and \(b\) is given by
 \[
  \ceff_{ab}=\frac{i_a-i'_a}{u_a-u'_a-u_b+u'_b}=\frac{D^k}{H^k_{ab}+H^k_{ba}}=\frac{D^k}{K^k_{ab}},
 \]
 which completes the proof.
\end{proof}

\subsection{A non-monotone example}

We have seen many nice properties of the network. The next step in the reversible case is making use of Rayleigh's monotonicity property: in the reversible case the effective resistance is a non-decreasing function of any of the individual resistances. Here we show an example to demonstrate that this is not the case in the irreversible case, the naive approach does not work. Resistance values below are in Ohms.
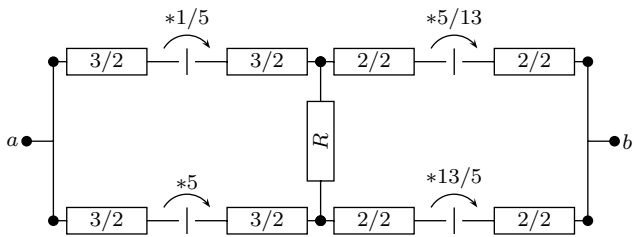
\begin{center}
 \begin{pspicture}(230,100)
  \qdisk(5,40){2}
  \psline(5,40)(15,40)(15,70)
  \psline(15,40)(15,10)
  \rput(15,70){\alap{3}{1/5}}
  \rput(115,70){\alap{2}{5/13}}
  \rput(15,10){\alap{3}{5}}
  \rput(115,10){\alap{2}{13/5}}
  \psline(215,70)(215,40)(225,40)
  \psline(215,10)(215,40)
  \qdisk(225,40){2}
  \rput{90}(115,10){\resist{R}}
  \uput{3}[180](5,40){\footnotesize\(a\)}
  \uput{3}[0](225,40){\footnotesize\(b\)}
 \end{pspicture}
\end{center}
We immediately rewrite this network to an equivalent form using the primer and secunder alternatives:
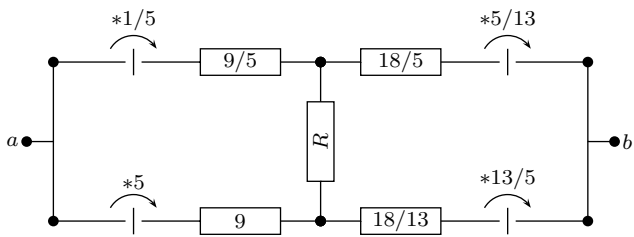
\begin{center}
 \begin{pspicture}(230,100)
  \qdisk(5,40){2}
  \psline(5,40)(15,40)(15,70)
  \psline(15,40)(15,10)
  \rput(15,70){\secunderno{9/5}{1/5}}
  \rput(115,70){\primerno{18/5}{5/13}}
  \rput(15,10){\secunderno{9}{5}}
  \rput(115,10){\primerno{18/13}{13/5}}
  \psline(215,70)(215,40)(225,40)
  \psline(215,10)(215,40)
  \qdisk(225,40){2}
  \rput{90}(115,10){\resist{R}}
  \uput{3}[180](5,40){\footnotesize\(a\)}
  \uput{3}[0](225,40){\footnotesize\(b\)}
 \end{pspicture}
\end{center}
First notice that a circular current of 4 Amperes in the positive direction and no current through \(R\) gives a constant 9 Volts free solution, thus the network is Markovian for all \(R\) values. Therefore it has an effective resistance, and it is perhaps easiest to compute if we fix \(u_a=5\) Volts and \(u_b=0\). Then we just need to figure out currents in
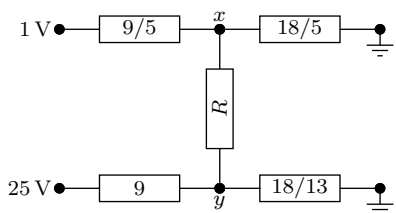
\begin{center}
 \begin{pspicture}(140,100)
  \qdisk(10,70){2}
  \qdisk(10,10){2}
  \qdisk(130,70){2}
  \qdisk(130,10){2}
  \rput(10,70){\resist{9/5}}
  \rput(10,10){\resist{9}}
  \rput(70,70){\resist{18/5}}
  \rput(70,10){\resist{18/13}}
  \rput{90}(70,10){\resist{R}}
  \psline(130,70)(130,64)
  \psline(125,64)(135,64)
  \psline(127,62)(133,62)
  \psline(129,60)(131,60)
  \psline(130,10)(130,4)
  \psline(125,4)(135,4)
  \psline(127,2)(133,2)
  \psline(129,0)(131,0)
  \uput{3}[180](10,70){\footnotesize{1\,V}}
  \uput{3}[180](10,10){\footnotesize{25\,V}}
  \uput{3}[90](70,70){\footnotesize\(x\)}
  \uput{3}[270](70,10){\footnotesize\(y\)}
 \end{pspicture}
\end{center}
One way of proceeding is to write the equations for the voltage dividers in \(x\) and \(y\). These are:
\[
 \ba
  u_x&=\frac{1\cdot\frac59+0\cdot\frac5{18}+u_y\cdot\frac1R}{\frac59+\frac5{18}+\frac1R}=\frac{10R}{15R+18}+\frac6{5R+6}\cdot u_y,\\
  u_y&=\frac{25\cdot\frac19+0\cdot\frac{13}{18}+u_x\cdot\frac1R}{\frac19+\frac{13}{18}+\frac1R}=\frac{50R}{15R+18}+\frac6{5R+6}\cdot u_x.
 \ea
\]
The solution is \(u_x=\frac{10R+72}{15R+36}\) and \(u_y=\frac{50R+72}{15R+36}\), and the effective resistance is
\[
 \reff_{ab}=\frac{u_a}{i_a}=\frac{u_a}{-i_b}=\frac{u_a}{\frac5{18}u_x+\frac{13}{18}u_y}=\frac{675R+1620}{350R+648}=\frac{27}{14}+\frac{1296}{1225R+2268},
\]
a decreasing function of \(R\). The situation reminds the authors to Braess's paradox \cite{braess}.

\subsection{The Dirichlet and Thomson principles}

Having the direct monotonicity approach failed, we now give an insight to the Dirichlet and Thomson energy minimum principles for the irreversible case. These are the fundamental principles that enable one to derive Rayleigh's monotonicity law in the reversible case. The irreversible case was established in this form by Gaudilli\`ere-Landim and Slowik \cite{gau-lan,slowik_var}, below we simply give a translation of their results without proof into the electrical language. In fact we stick to the notation of Slowik's Proposition 2.6 as closely as possible.

As before, take two sets \(A\cap B=\emptyset\) of vertices, and define
\[
 \mathcal H_{A,B}:\,=\{u:V\to\Rb\,:\,u|_A\equiv1,\ u|_B\equiv0\},
\]
we think of such functions as voltages with respective boundary conditions on \(A\) and \(B\). Set \(\mathfrak U^0_{A,B}\) as the set of currents \(i_{\cdot\cdot}\) with zero external currents \(i_x\) for \(x\notin(A\cup B)\) (c.f.\ \eqref{eq:ic}), and total incoming current \(i_A=0=i_B\) to the set \(A\) (and therefore to the set \(B\) as well), see \eqref{eq:tic}. The next quantity to define is, for currents \(i_{\cdot\cdot}\),
\[
 \mathcal D(i)=\frac12\sum_{y\sim x\in V}\frac1{\mu_xP^\text s_{xy}}i^2_{xy}=\frac12\sum_{y\sim x\in V}R_{xy}i^2_{xy},
\]
the ohmic power losses on the resistors, see \eqref{eq:ps}. Finally, recall \eqref{eq:iu}, and reverse the amplifiers in there to get \(\hat i^\cdot\). We can now state
\begin{pr}[Dirichlet principle, Slowik \cite{slowik_var} Proposition 2.6]
 \[
  \cp(A,\,B)=\min_{u\in\mathcal H_{A,B}}\min_{i\in\mathfrak U^0_{A,B}}\mathcal D(\hat i^u-i).
 \]
 The minimum is attained for \(u=\frac12(u_{AB}+\hat u_{AB})\) and \(i=\hat i^u-i^{\text su_{AB}}\), where \(u_{AB}\) and \(\hat u_{AB}\) are the physical potentials in the network and in the reversed network under our boundary conditions, respectively, and \(i^{\text su_{AB}}\) is the current that would result under the potential \((u_{AB})_\cdot\) without the amplifiers.
\end{pr}
In words: find a potential function \(u\) with our boundary conditions (this results in currents \(\hat i^u\) in the reversed network) and a divergence free current \(i\) on the free vertices with total incoming flow \(i_A=0\) such that the difference of these two currents minimises the ohmic losses on the network. Then these ohmic losses sum up to the total physical power required to maintain the boundary conditions (this is the effective conductance \(\ceff_{AB}=\cp(A,\,B)\) since the boundary voltage difference \(U_A-U_B=1\)). We emphasise that the minimisers \(u\) and \(i\) are non-physical, except in the reversible case when \(u\) becomes the physical voltage while \(i\equiv0\).

Next, we define
\[
 G^0_{A,B}:\,=\{u:V\to\Rb\,:\,u|_A\equiv u|_B\equiv0\},
\]
and \(\mathfrak U^1_{A,B}\) as the set of currents \(i_{\cdot\cdot}\) with zero external currents \(i_x\) for \(x\notin(A\cup B)\) (c.f.\ \eqref{eq:ic}), and total incoming current \(i_A=1=-i_B\) to the set \(A\) (and therefore -1 to the set \(B\)), see \eqref{eq:tic}.
\begin{pr}[Thomson principle, Slowik \cite{slowik_var} Proposition 2.6]
 \[
  \cp(A,\,B)=\max_{i\in\mathfrak U^1_{A,B}}\max_{u\in G^0_{A,B}}\frac1{\mathcal D(i-\hat i^u)}.
 \]
 The maximum is attained for \(u=\frac12(u_{AB}-\hat u_{AB})/\cp(A,\,B)\), and \(i=i_0^{\text su_{AB}}+\hat i^u\), where \(i_0^{\text su_{AB}}\) is the current that would result under the potential \((u_{AB})_\cdot\) without the amplifiers, except that it is normalised to have unit total inflow in \(A\).
\end{pr}
In words: find a potential function \(u\) that vanishes in \(A\) and \(B\) (this results in currents \(\hat i^u\) in the reversed network), and a unit flow \(i\) such that the difference of these currents minimises the ohmic losses on the network. Then these ohmic losses sum up to the total physical power required to maintain a unit flow (the reciprocal is again the effective conductance since the total current flow is one). Again, the minimisers are non-physical, except for the reversible case when \(u\equiv0\) and \(i\) is the physical unit current flow.

How these principles can be used towards monotonicity is a question left for future work.

\section{Series, parallel and star-delta transformations}

In this last section we dive a little bit into ``network-algebra'' by showing how units in series, parallel, star or delta configurations behave. Series and parallel are simple and nice situations:
\begin{pr}
 Two of our units in series of respective parameters \((R,\,\la)\) and \((Q,\,\mu)\) can be replaced by a single unit of parameters
 \[
  \Bigl(R\frac{(\lambda+1)\mu}{\lambda\mu+1}+Q\frac{\mu+1}{\lambda\mu+1},\ \la\mu\Bigr).
 \]
\end{pr}
\begin{proof}
 We use the primer and secunder alternatives as
 
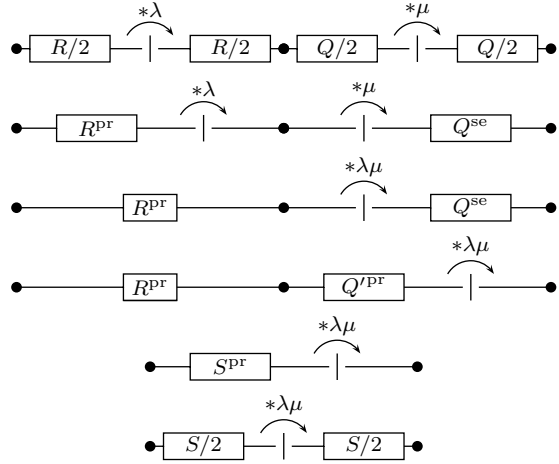
\begin{center}
  \begin{pspicture}(200,190)
   \rput(0,160){\alap{R}{\la}}
   \rput(100,160){\alap{Q}{\mu}}
 
   \rput(0,130){\primer{R}{\la}}
   \rput(100,130){\secunder{Q}{\mu}}
 
   \qdisk(0,100){2}
   \psline(0,100)(40,100)(40,105)(60,105)(60,95)(40,95)(40,100)
   \rput(50,100){\footnotesize\(R^\text{pr}\)}
   \psline(60,100)(100,100)
   \rput(100,100){\secunder{Q}{\la\mu}}
 
   \qdisk(0,70){2}
   \psline(0,70)(40,70)(40,75)(60,75)(60,65)(40,65)(40,70)
   \psline(40,65)(60,65)
   \rput(50,70){\footnotesize\(R^\text{pr}\)}
   \psline(60,70)(100,70)
   \rput(100,70){\primer{{Q'}}{\la\mu}}
   
   \rput(50,40){\primer{S}{\la\mu}}
 
   \rput(50,10){\alap{S}{\la\mu}}
  \end{pspicture}
 \end{center}
 It is obvious that the parameter of the voltage amplifier is \(\la\mu\). Applying the transcription formula for each step, the resistance of the substitute element can be determined:
 \[
  \ba
   S=S^\text{pr}\frac{2\lambda\mu}{\lambda\mu+1}=(R^\text{pr}+{Q'}^\text{pr})\frac{2\lambda\mu}{\lambda\mu+1}&=R^\text{pr}\frac{2\lambda\mu}{\lambda\mu+1}+Q^\text{se}\frac{2}{\lambda\mu+1}\\
   &=R\frac{(\lambda+1)\mu}{\lambda\mu+1}+Q\frac{\mu+1}{\lambda\mu+1}.
  \ea
 \]
\end{proof}
\begin{pr}
 Two of our units in parallel of respective parameters \((R,\,\la)\) and \((Q,\,\mu)\) can be replaced by a single unit of parameters
 \[
  \Bigl(\frac{RQ}{R+Q},\ \frac{Q(\mu+1)}{Q(\mu+1)+R(\la+1)}\cdot\la+\frac{R(\la+1)}{Q(\mu+1)+R(\la+1)}\cdot\mu\Bigr).
 \]
\end{pr}
Notice the classical parallel formula for the resistance, and the weighted average for the amplifier.
\begin{proof}
 This case cannot be reduced with transformations into one single unit, but the alternative elements are still useful. Below are two equivalent circuits.
 
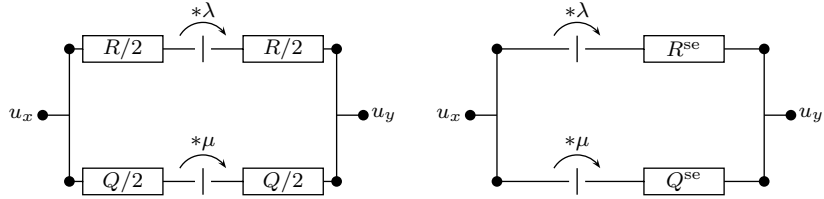
\begin{center}
  \begin{pspicture}(300,110)
   \rput(20,30){\alap{Q}{\mu}}
   \rput(20,80){\alap{R}{\lambda}}
   \qdisk(10,55){2}
   \qdisk(130,55){2}
   \psline(20,30)(20,80)
   \psline(120,30)(120,80)
   \psline(10,55)(20,55)
   \psline(120,55)(130,55)
   \uput{3}[180](10,55){\footnotesize\(u_x\)}
   \uput{3}[0](130,55){\footnotesize\(u_y\)}
 
   \rput(180,30){\secunder{Q}{\mu}}
   \rput(180,80){\secunder{R}{\lambda}}
   \qdisk(170,55){2}
   \qdisk(290,55){2}
   \psline(180,30)(180,80)
   \psline(280,30)(280,80)
   \psline(170,55)(180,55)
   \psline(280,55)(290,55)
   \uput{3}[180](170,55){\footnotesize\(u_x\)}
   \uput{3}[0](290,55){\footnotesize\(u_y\)}
  \end{pspicture}
 \end{center}
 The total current from \(x\) to \(y\) will be the sum of currents like in \eqref{eq:arram} for the top and the bottom branches, therefore it will be in the same form. This proves that a single unit can be used as a replacement. Its secunder alternative will be
 
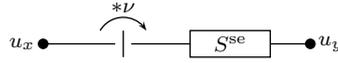
\begin{center}
  \begin{pspicture}(120,50)
   \rput(10,20){\secunder{S}{\nu}}
   \uput{3}[180](10,20){\footnotesize\(u_x\)}
   \uput{3}[0](110,20){\footnotesize\(u_y\)}
  \end{pspicture}
 \end{center}
 It remains to determine the parameters \(S\) and \(\nu\). Assume first \(u_x=1\) and zero total current that is, leave vertex \(y\) free. Then the secunder resistors act as a voltage divider, giving
 \[
  u_y=\frac{\la Q^\text{se}+\mu R^\text{se}}{Q^\text{se}+R^\text{se}}.
 \]
 In the simple unit this agrees to the value \(\nu\) of the amplifier, therefore
 \[
  \nu=\frac{Q^\text{se}\la+R^\text{se}\mu}{Q^\text{se}+R^\text{se}}=\frac{Q\la(\mu+1)+R\mu(\la+1)}{Q(\mu+1)+R(\la+1)}.
 \]
 Next, when \(u_x=0\), the amplifiers keep the potentials at zero, and the parallel formula
 \[
  S^\text{se}=\frac{R^\text{se}Q^\text{se}}{R^\text{se}+Q^\text{se}}
 \]
 follows. Returning to the original alternative,
 \[
  S=\frac2{\nu+1}\cdot S^\text{se}=2\frac{Q(\mu+1)+R(\la+1)}{(R+Q)(\la+1)(\mu+1)}\cdot\frac{R^\text{se}Q^\text{se}}{R^\text{se}+Q^\text{se}}=\frac{RQ}{R+Q}.
 \]
\end{proof}
Notice that in both the series or parallel formulas the resulting resistances are monotone increasing functions of the original ones. Not all networks can, however, be reduced using only series or parallel substitutions. The next step of transformations for classical resistor networks is the star-delta transformation. As we will see shortly it is here where non-monotonicity issues begin. Our non-monotone example is also one that cannot be reduced using only series or parallel substitutions.

In our case, star and delta look like
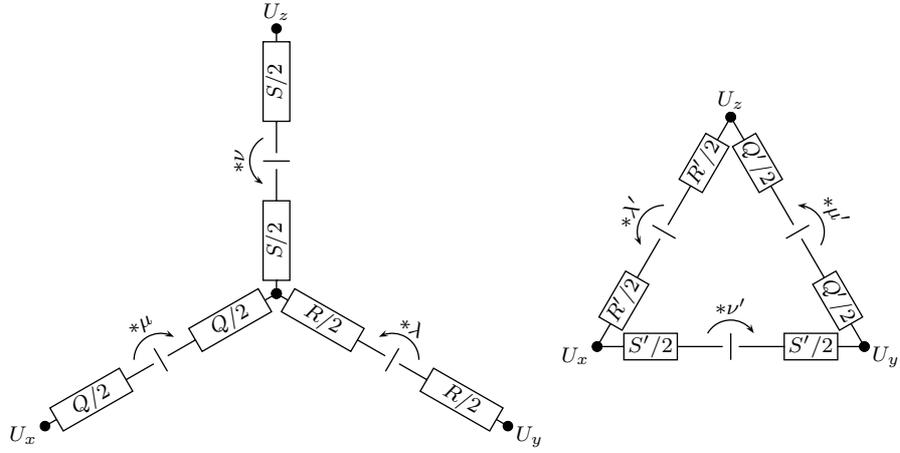
\begin{center}
 \begin{pspicture}(340,170)
  \rput{90}(110,60){\alapf{S}{\nu}}
  \rput{330}(110,60){\alapf{R}{\la}}
  \rput{30}(23.4,10){\alap{Q}{\mu}}
  \uput{3}[210](23.4,10){\footnotesize\(U_x\)}
  \uput{3}[330](196.6,10){\footnotesize\(U_y\)}
  \uput{3}[90](110,160){\footnotesize\(U_z\)}

  \rput(230,40){\dalap{S'}{\nu'}}
  \rput{60}(230,40){\dalapf{R'}{\la'}}
  \rput{300}(280,126.6){\dalapf{Q'}{\mu'}}
  \uput{3}[210](230,40){\footnotesize\(U_x\)}
  \uput{3}[330](330,40){\footnotesize\(U_y\)}
  \uput{3}[90](280,126.6){\footnotesize\(U_z\)}
 \end{pspicture}
\end{center}
where it is essential that the centre of star has no further connections. The question is whether the parameters can be linked so that these two networks behave identically under all scenarios. We start by rewriting the above into the equivalent secunder alternatives, and work with those thereafter. Any formulas can be rewritten into the original parameters via \eqref{eq:prisech}, we avoid that for the sake of simplicity (well...).
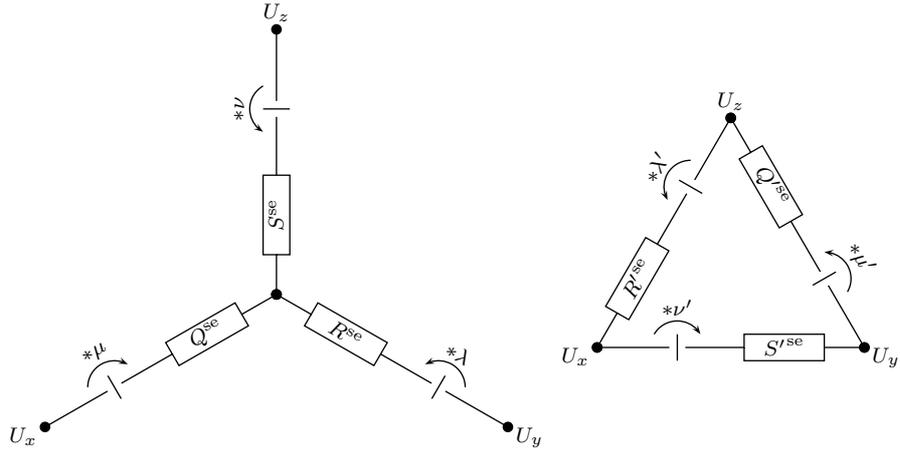
\begin{center}
 \begin{pspicture}(340,170)
  \rput{90}(110,60){\secunderf{S}{\nu}}
  \rput{330}(110,60){\secunderf{R}{\la}}
  \rput{30}(23.4,10){\secunder{Q}{\mu}}
  \uput{3}[210](23.4,10){\footnotesize\(U_x\)}
  \uput{3}[330](196.6,10){\footnotesize\(U_y\)}
  \uput{3}[90](110,160){\footnotesize\(U_z\)}

  \rput(230,40){\secunder{{S'}}{\nu'}}
  \rput{60}(230,40){\secunderf{{R'}}{\la'}}
  \rput{300}(280,126.6){\secunderf{{Q'}}{\mu'}}
  \uput{3}[210](230,40){\footnotesize\(U_x\)}
  \uput{3}[330](330,40){\footnotesize\(U_y\)}
  \uput{3}[90](280,126.6){\footnotesize\(U_z\)}
 \end{pspicture}
\end{center}
With the notations of this picture, we have
\begin{pr}\label{pr:std}
 Any star can be transformed into an equivalent delta, the parameters of which are given by
 \[
  \ba
   {S'}^\text{se}&=\frac{R^\text{se}S^\text{se}+Q^\text{se}S^\text{se}+Q^\text{se}R^\text{se}}{\la S^\text{se}},\\
   {Q'}^\text{se}&=\frac{R^\text{se}S^\text{se}+Q^\text{se}S^\text{se}+Q^\text{se}R^\text{se}}{\nu Q^\text{se}},\\
   {R'}^\text{se}&=\frac{R^\text{se}S^\text{se}+Q^\text{se}S^\text{se}+Q^\text{se}R^\text{se}}{\mu R^\text{se}},
  \ea
 \]
 and \(\la'=\frac\nu\mu\), \(\nu'=\frac\mu\la\), \(\mu'=\frac\la\nu\).
\end{pr}
Not every delta, however, can be transformed into a star.
\begin{pr}
 A delta can be transformed into an equivalent star if and only if 
 \be
  \la'\cdot\nu'\cdot\mu'=1.\label{eq:dprod}
 \ee
 Even in this case the resulting star is not unique: with any positive number \(\al>0\), it can have parameters
 \[
  \ba
   S^\text{se}&=\frac{{\la'}^{1/3}{\nu'}^{2/3}\al{R'}^\text{se}{Q'}^\text{se}}{{\la'}^{2/3}{\nu'}^{1/3}{Q'}^\text{se}+{\nu'}^{2/3}{\mu'}^{1/3}{R'}^\text{se}+{\la'}^{1/3}{\mu'}^{2/3}{S'}^\text{se}},\\
   Q^\text{se}&=\frac{{\nu'}^{1/3}{\mu'}^{2/3}\al{R'}^\text{se}{S'}^\text{se}}{{\la'}^{2/3}{\nu'}^{1/3}{Q'}^\text{se}+{\nu'}^{2/3}{\mu'}^{1/3}{R'}^\text{se}+{\la'}^{1/3}{\mu'}^{2/3}{S'}^\text{se}},\\
   R^\text{se}&=\frac{{\la'}^{2/3}{\mu'}^{1/3}\al{Q'}^\text{se}{S'}^\text{se}}{{\la'}^{2/3}{\nu'}^{1/3}{Q'}^\text{se}+{\nu'}^{2/3}{\mu'}^{1/3}{R'}^\text{se}+{\la'}^{1/3}{\mu'}^{2/3}{S'}^\text{se}},
  \ea
 \]
 and \(\la={\la'}^{1/3}{\mu'}^{2/3}\al\), \(\nu={\la'}^{2/3}{\nu'}^{1/3}\al\), \(\mu={\nu'}^{2/3}{\mu'}^{1/3}\al\).
\end{pr}
\begin{proof}[Proof of both star\(\to\)delta and delta\(\to\)star]
 We determine and compare the incoming currents on vertices \(x,\,y\) and \(z\) in the two networks. We start with star. The voltages at the outer points of the resistances are \(\mu U_x\), \(\la U_y\) and \(\nu U_z\), thus the voltage in the centre point is
 \[
 U=\frac{\mu U_x\frac1{Q^\text{se}}+\la U_y\frac1{R^\text{se}}+\nu U_z\frac1{S^\text{se}}}{\frac1{Q^\text{se}}+\frac1{R^\text{se}}+\frac1{S^\text{se}}}=\frac{\mu U_xR^\text{se}S^\text{se}+\la U_yQ^\text{se}S^\text{se}+\nu U_zQ^\text{se}R^\text{se}}{R^\text{se}S^\text{se}+Q^\text{se}S^\text{se}+Q^\text{se}R^\text{se}}.
 \]
 Therefore, the respective currents flowing from \(x\), \(y\) and \(z\) into the centre point are
 \[
  \ba
   i_x&=\frac{\mu U_x-U}{Q^\text{se}}=\frac{\mu(S^\text{se}+R^\text{se})U_x-\la S^\text{se}U_y-\nu R^\text{se}U_z}{R^\text{se}S^\text{se}+Q^\text{se}S^\text{se}+Q^\text{se}R^\text{se}},\\
   i_y&=\frac{\la U_y-U}{R^\text{se}}=\frac{\la(S^\text{se}+Q^\text{se})U_y-\mu S^\text{se}U_x-\nu Q^\text{se}U_z}{R^\text{se}S^\text{se}+Q^\text{se}S^\text{se}+Q^\text{se}R^\text{se}},\\
   i_z&=\frac{\nu U_z-U}{S^\text{se}}=\frac{\nu (R^\text{se}+Q^\text{se})U_z-\mu R^\text{se}U_x-\la Q^\text{se}U_y}{R^\text{se}S^\text{se}+Q^\text{se}S^\text{se}+Q^\text{se}R^\text{se}}.
  \ea
 \]
 Next we turn to delta. The currents flowing on the edges are:
 \[
  i_{xy}=\frac{\nu'U_x-U_y}{{S'}^\text{se}},\qquad i_{yz}=\frac{\mu'U_y-U_z}{{Q'}^\text{se}},\qquad i_{zx}=\frac{\la'U_z-U_x}{{R'}^\text{se}}.
 \]
 Thus, the currents flowing from $x$, $y$ and $z$ into the network can be written as
 \[
  \ba
   i_x&=i_{xy}-i_{zx}=\frac{(\nu'{R'}^\text{se}+{S'}^\text{se})U_x-{R'}^\text{se}U_y-\la'{S'}^\text{se}U_z}{{R'}^\text{se}{S'}^\text{se}},\\
   i_y&=i_{yz}-i_{xy}=\frac{(\mu'{S'}^\text{se}+{Q'}^\text{se})U_y-\nu'{Q'}^\text{se}U_x-{S'}^\text{se}U_z}{{Q'}^\text{se}{S'}^\text{se}},\\
   i_z&=i_{zx}-i_{yz}=\frac{(\la'{Q'}^\text{se}+{R'}^\text{se})U_z-{Q'}^\text{se}U_x-\mu'{R'}^\text{se}U_y}{{Q'}^\text{se}{R'}^\text{se}}.
  \ea
 \]
 In  a star\,/\,delta substitution the currents have to equal for all possible voltages $U_x$, $U_y$ and $U_z$. Hence, by comparing the coefficients of the voltages in the formulas for the currents, the connections between the quantities can be determined. It is subservient to consider first the coefficient of $U_y$ in the formula for $i_x$, the coefficient of $U_z$ in the formula for $i_y$ and the coefficient of $U_x$ in the formula for $i_z$:
 \be
  \ba
   \frac1{{S'}^\text{se}}&=\frac{\la S^\text{se}}{R^\text{se}S^\text{se}+Q^\text{se}S^\text{se}+Q^\text{se}R^\text{se}},\\
   \frac1{{Q'}^\text{se}}&=\frac{\nu Q^\text{se}}{R^\text{se}S^\text{se}+Q^\text{se}S^\text{se}+Q^\text{se}R^\text{se}},\\
   \frac1{{R'}^\text{se}}&=\frac{\mu R^\text{se}}{R^\text{se}S^\text{se}+Q^\text{se}S^\text{se}+Q^\text{se}R^\text{se}}.
  \ea\label{eq:ellcsd}
 \ee
 Second, consider the coefficient of $U_z$ in $i_x$, the coefficient of $U_x$ in $i_y$ and the coefficient of $U_y$ in $i_z$:
 \[
  \ba
   \frac{\la'}{{R'}^\text{se}}&=\frac{\nu R^\text{se}}{R^\text{se}S^\text{se}+Q^\text{se}S^\text{se}+Q^\text{se}R^\text{se}},\\
   \frac{\nu'}{{S'}^\text{se}}&=\frac{\mu S^\text{se}}{R^\text{se}S^\text{se}+Q^\text{se}S^\text{se}+Q^\text{se}R^\text{se}},\\
   \frac{\mu'}{{Q'}^\text{se}}&=\frac{\la Q^\text{se}}{R^\text{se}S^\text{se}+Q^\text{se}S^\text{se}+Q^\text{se}R^\text{se}}.
  \ea
 \]
 By dividing the corresponding equations in the two triplets of equations we get:
 \[
  \la'=\frac\nu\mu,\qquad\nu'=\frac\mu\la\text{ and }\qquad\mu'=\frac\la\nu.
 \]
 Finally, with these substitutions the coefficients of \(U_x\) in \(i_x\), \(U_y\) in \(i_y\) and \(U_z\) in \(i_z\) also match. This proves Proposition \ref{pr:std}.
 
 Notice that for any star, in the substitute delta we have \(\la'\cdot\nu'\cdot\mu'=1\) and that multiplying the amplifiers in the star by a common constant does not change the parameters of the voltage amplifiers in the substitute delta. Therefore, the inversion of the previous formulas is possible only if \eqref{eq:dprod} holds, and even in this case it is not uniquely determined, a constant factor has to be chosen. This can be written as \(\la\nu\mu=\alpha^3\) where $\alpha>0$ is the free parameter. Then
 \[
  \la={\la'}^{1/3}{\mu'}^{2/3}\al,\qquad\nu={\la'}^{2/3}{\nu'}^{1/3}\al\text{ and }\qquad\mu={\nu'}^{2/3}{\mu'}^{1/3}\al.
 \]
 To invert the resistances \eqref{eq:ellcsd} first note that
 \[
  \frac1{\mu{R'}^\text{se}}\cdot\frac1{\la{S'}^\text{se}}+\frac1{\nu{Q'}^\text{se}}\cdot\frac1{\la{S'}^\text{se}}+\frac1{\nu{Q'}^\text{se}}\cdot\frac1{\mu{R'}^\text{se}}=\frac1{R^\text{se}S^\text{se}+Q^\text{se}S^\text{se}+Q^\text{se}R^\text{se}}.
 \]
 Thus from \eqref{eq:ellcsd} and from \eqref{eq:dprod}:
 \[
  \ba
   S^\text{se}&=\frac{\frac1{\la{S'}^\text{se}}}{\frac1{\mu{R'}^\text{se}}\cdot\frac1{\la{S'}^\text{se}}+\frac1{\nu{Q'}^\text{se}}\cdot\frac1{\la{S'}^\text{se}}+\frac1{\nu{Q'}^\text{se}}\cdot\frac1{\mu{R'}^\text{se}}}\\
   &=\frac{{\la'}^{1/3}{\nu'}^{2/3}\al{R'}^\text{se}{Q'}^\text{se}}{{\la'}^{2/3}{\nu'}^{1/3}{Q'}^\text{se}+{\nu'}^{2/3}{\mu'}^{1/3}{R'}^\text{se}+{\la'}^{1/3}{\mu'}^{2/3}{S'}^\text{se}},\\
   Q^\text{se}&=\frac{\frac1{\nu{Q'}^\text{se}}}{\frac1{\mu{R'}^\text{se}}\cdot\frac1{\la{S'}^\text{se}}+\frac1{\nu{Q'}^\text{se}}\cdot\frac1{\la{S'}^\text{se}}+\frac1{\nu{Q'}^\text{se}}\cdot\frac1{\mu{R'}^\text{se}}}\\
   &=\frac{{\nu'}^{1/3}{\mu'}^{2/3}\al{R'}^\text{se}{S'}^\text{se}}{{\la'}^{2/3}{\nu'}^{1/3}{Q'}^\text{se}+{\nu'}^{2/3}{\mu'}^{1/3}{R'}^\text{se}+{\la'}^{1/3}{\mu'}^{2/3}{S'}^\text{se}},\\
   R^\text{se}&=\frac{\frac1{\mu{R'}^\text{se}}}{\frac1{\mu{R'}^\text{se}}\cdot\frac1{\la{S'}^\text{se}}+\frac1{\nu{Q'}^\text{se}}\cdot\frac1{\la{S'}^\text{se}}+\frac1{\nu{Q'}^\text{se}}\cdot\frac1{\mu{R'}^\text{se}}}\\
   &=\frac{{\la'}^{2/3}{\mu'}^{1/3}\al{Q'}^\text{se}{S'}^\text{se}}{{\la'}^{2/3}{\nu'}^{1/3}{Q'}^\text{se}+{\nu'}^{2/3}{\mu'}^{1/3}{R'}^\text{se}+{\la'}^{1/3}{\mu'}^{2/3}{S'}^\text{se}}.
  \ea
 \]
\end{proof}
The condition \eqref{eq:dprod} is that at any constant potential the delta has no circular current by itself. This is rather restrictive, thus delta\(\to\)star transformations cannot be used to reduce a general network. After the lack of monotonicity, this is the second serious drawback of our networks compared to the classical resistor-only case.

\section*{Acknowledgements}

The authors wish to thank Edward Crane, Nic Freeman, Alexandre Gaudilli\`ere, Claudio Landim, G\'abor Pete, Martin Slowik, Andr\'as Telcs and B\'alint T\'oth for stimulating discussions on this project.

\bibliographystyle{plain}
\bibliography{refsmarton}
\end{document}